\documentclass[12pt,thmsa]{article}
\usepackage{amssymb}
\usepackage{amsfonts}
\usepackage{sw20bams}

\input tcilatex
\begin{document}

\author{Steven R. Finch}
\title{Planar Projections and Second Intrinsic Volume}
\date{March 12, 2012}
\maketitle

\begin{abstract}
Consider random shadows of a cube and of a regular tetrahedron. Area and
perimeter of the former are positively dependent (with correlation $0.915...$%
), whereas area and perimeter of the latter appear to be negatively
dependent. \ This is only one result of many, all involving generalizations
of mean width.
\end{abstract}

\footnotetext{%
Copyright \copyright\ 2012 by Steven R. Finch. All rights reserved.}Let $C$
be a convex body in $\mathbb{R}^{n}$. Let $S$ be a $1$-dimensional subspace
passing through the origin. A \textbf{width} is the length of the orthogonal
projection of $C$ on $S$. If $S$ is uniformly distributed on the
Grassmannian manifold $\mathbb{G}^{n,1}$, then the width $w(S)$ is a random
variable. \ First and second moments of $w$ are known for $C=$ the regular $%
n $-simplex, $C=$ the $n$-cube and $C=$ the regular $n$-crosspolytope, each
centered at the origin \cite{Fi1, Fi2}. The \textbf{mean width} is sometimes
called the \textbf{mean linear projection} \cite{Ms}-- a slight abuse of
language -- because a projection of $C$ on $S$ has only one parameter of
interest: its length.

Let instead $S$ be a $2$-dimensional subspace passing through the origin.
Orthogonal projections of $C$ on $S$ now have three parameters of interest.
\ One parameter is the number of polygonal vertices; associated
probabilities will be mentioned later. The other two parameters constitute
direct generalizations of width and will be our main focus:

\begin{itemize}
\item a \textbf{chorowidth} is the polygonal area of the orthogonal
projection of $C$ on $S$ (\textquotedblleft choro\textquotedblright\ is
Greek for place or area);

\item a \textbf{periwidth} is the polygonal circumference of the orthogonal
projection of $C$ on $S$ (\textquotedblleft peri\textquotedblright\ is Greek
for around or enclosing).
\end{itemize}

\noindent If $S$ is uniformly distributed on $\mathbb{G}^{n,2}$, then
chorowidth $cw(S)$ and periwidth $pw(S)$ are random variables. We shall
examine joint moments of $(cw,pw)$ for the same regular polytopes as before,
for the special cases $n=3$ and $n=4$.

\section{Method for Computing Moments}

To generate a random $2$-subspace $S$ in $\mathbb{R}^{3}$, we select a
random point $U$ uniformly on the $2$-sphere of unit radius. The desired
subspace is the set of all vectors orthogonal to $U$.

We then project the (fixed) convex polyhedron $C$ orthogonally onto $S$.
This is done by forming the convex hull of images of all vertices of $C$. If 
$C=$ the regular $3$-simplex, the resultant polygon in the plane has $3$ or $%
4$ vertices almost surely. If $C=$ the $3$-cube, there are $6$ vertices
almost surely; if $C=$ the regular $3$-crosspolytope, there are $4$ or $6$
vertices almost surely. We compute the area and circumference of this
polygon via standard formulas.

More precisely, if $U=(x,y,z)$ is of unit length, then the matrix%
\[
M_{3}=\left( 
\begin{array}{ccc}
\sqrt{1-x^{2}} & -\frac{x\,y}{\sqrt{1-x^{2}}} & -\frac{x\,z}{\sqrt{1-x^{2}}}
\\ 
0 & \frac{z}{\sqrt{1-x^{2}}} & -\frac{y}{\sqrt{1-x^{2}}} \\ 
0 & 0 & 0%
\end{array}%
\right) 
\]%
projects $C$ orthogonally onto a plane, rotated in $\mathbb{R}^{3}$ to
coincide with the $2$-subspace spanned by $(1,0,0)$ and $(0,1,0)$ for
convenience. \ Let $T$ be the $3$-row matrix whose columns constitute all
vertices of $C$. Then the first $2$ rows of $M_{3}T$ constitute all images
of the vertices in $\mathbb{R}^{2}$ and $2$-dimensional convex hull
algorithms apply naturally.

Spherical coordinates in $\mathbb{R}^{3}$: 
\[
\begin{array}{ccccc}
x=\cos \theta \sin \varphi , &  & y=\sin \theta \sin \varphi , &  & z=\cos
\varphi%
\end{array}%
\]%
will be used throughout for $U$, where $0\leq \theta <2\pi $, $0\leq \varphi
\leq \pi $. \ The corresponding Jacobian determinant is $\sin \varphi $; it
is best to think of $(\theta ,\varphi )$ as possessing joint density $\frac{1%
}{4\pi }\sin \varphi $.

Moving up a dimension, to generate a random $2$-subspace $S$ in $\mathbb{R}%
^{4}$, we first select a random point $U$ uniformly on the $3$-sphere of
unit radius. The set of unit vectors orthogonal to $U$ form a $2$-sphere; we
select a random point $V$ uniformly on this $2$-sphere. The desired subspace
is then the set of all vectors orthogonal to both $U$ and $V$.

We next project the (fixed) convex polyhedron $C$ orthogonally onto $S$.
This is done by forming the convex hull of images of all vertices of $C$. If 
$C=$ the regular $4$-simplex, the resultant polygon in the plane has $3$, $4$
or $5$ vertices almost surely. If $C=$ the $4$-cube, there are $8$ vertices
almost surely; if $C=$ the regular $4$-crosspolytope, there are $4$, $6$ or $%
8$ vertices almost surely. Area and circumference of this polygon are
computed as before.

More precisely, if $U=(x,y,z,w)$ and $V=(p,q,r,s)$ are orthogonal and of
unit length, then the matrix%
\[
M_{4}=\left( 
\begin{array}{cccc}
\sqrt{1-p^{2}-x^{2}} & -\frac{p\,q+x\,y}{\sqrt{1-p^{2}-x^{2}}} & -\frac{%
p\,r+x\,z}{\sqrt{1-p^{2}-x^{2}}} & -\frac{p\,s+x\,w}{\sqrt{1-p^{2}-x^{2}}}
\\ 
0 & \frac{r\,w-s\,z}{\sqrt{1-p^{2}-x^{2}}} & -\frac{q\,w-s\,y}{\sqrt{%
1-p^{2}-x^{2}}} & \frac{q\,z-r\,y}{\sqrt{1-p^{2}-x^{2}}} \\ 
0 & 0 & 0 & 0 \\ 
0 & 0 & 0 & 0%
\end{array}%
\right) 
\]%
projects $C$ orthogonally onto a plane, rotated in $\mathbb{R}^{4}$ to
coincide with the $2$-subspace spanned by $(1,0,0,0)$ and $(0,1,0,0)$ for
convenience. \ Let $T$ be the $4$-row matrix whose columns constitute all
vertices of $C$. Then the first $2$ rows of $M_{4}T$ constitute all images
of the vertices in $\mathbb{R}^{2}$ and $2$-dimensional convex hull
algorithms again apply naturally.

Spherical coordinates in $\mathbb{R}^{4}$:%
\[
\begin{array}{ccccccc}
x=\cos \theta \sin \varphi \sin \psi , &  & y=\sin \theta \sin \varphi \sin
\psi , &  & z=\cos \varphi \sin \psi , &  & w=\cos \psi%
\end{array}%
\]%
will be used throughout for $U$, where $0\leq \theta <2\pi $, $0\leq \varphi
\leq \pi $, $0\leq \psi \leq \pi $. \ It follows that $V$ is given by%
\[
\left( 
\begin{array}{c}
p \\ 
q \\ 
r \\ 
s%
\end{array}%
\right) =\cos \kappa \sin \lambda \left( 
\begin{array}{c}
-y \\ 
x \\ 
-w \\ 
z%
\end{array}%
\right) +\sin \kappa \sin \lambda \left( 
\begin{array}{c}
-z \\ 
w \\ 
x \\ 
-y%
\end{array}%
\right) +\cos \lambda \left( 
\begin{array}{c}
-w \\ 
-z \\ 
y \\ 
x%
\end{array}%
\right) 
\]%
where $0\leq \kappa <2\pi $, $0\leq \lambda \leq \pi $. It is best to think
of $(\theta ,\varphi ,\psi ,\kappa ,\lambda )$ as possessing joint density $%
\frac{1}{2\pi ^{2}}\sin \varphi \sin ^{2}\psi \frac{1}{4\pi }\sin \lambda $.

\section{Three-Dimensional Results}

From the foregoing, we calculate mean chorowidth and mean periwidth via
integration:%
\[
\mathbb{E}\left( cw\right) =\dint\limits_{0}^{2\pi }\dint\limits_{0}^{\pi }%
\text{area}\left( \text{convex\_hull}\left( M_{3}T\right) \right) \frac{1}{%
4\pi }\sin \varphi \,d\varphi \,d\theta , 
\]%
\[
\mathbb{E}\left( pw\right) =\dint\limits_{0}^{2\pi }\dint\limits_{0}^{\pi }%
\text{circumference}\left( \text{convex\_hull}\left( M_{3}T\right) \right) 
\frac{1}{4\pi }\sin \varphi \,d\varphi \,d\theta 
\]%
and likewise for mean square chorowidth, mean square periwidth and joint
moment.

\subsection{3-Simplex (Tetrahedron)}

\[
\mathbb{E}\left( cw\right) =\frac{\sqrt{3}}{4}=0.433012701892219...=\frac{%
\text{surface area}}{4}, 
\]

\[
\mathbb{E}\left( cw^{2}\right) =\frac{1}{8}+\frac{\sqrt{2}}{4\pi }-\frac{1}{%
8\pi }\func{arcsec}(3)=0.188561220515812..., 
\]

\[
\mathbb{E}\left( pw\right) =\frac{3}{2}\left( \pi -\func{arcsec}(3)\right)
=2.865949854373527..., 
\]

\[
\mathbb{E}\left( pw^{2}\right) =8.2170808733..., 
\]%
\[
\mathbb{E}\left( cw\cdot pw\right) =1.2406348222... 
\]%
which imply that the correlation between $cw$ and $pw$ is $\approx -0.188$.
Closed-form expressions for $\mathbb{E}\left( pw^{2}\right) $ and $\mathbb{E}%
\left( cw\cdot pw\right) $ are not known. \ The fact that 
\[
\left( 
\begin{array}{c}
\text{surface area of a} \\ 
\text{3D convex body}%
\end{array}%
\right) =4\left( 
\begin{array}{c}
\text{mean area of} \\ 
\text{its 2D shadow}%
\end{array}%
\right) 
\]%
was first noticed by Cauchy \cite{SW, Chy1, Chy2, Chy3}. We also have \cite%
{AS, FK}%
\[
\mathbb{P}\left( \text{projection has 3 vertices}\right) =\frac{2}{\pi }%
\left( 3\func{arcsec}\left( 3\right) -\pi \right) =0.3509593121..., 
\]%
\[
\mathbb{P}\left( \text{projection has 4 vertices}\right) =\frac{3}{\pi }%
\left( \pi -2\func{arcsec}\left( 3\right) \right) =0.6490406878.... 
\]

\subsection{3-Cube}

\[
\mathbb{E}\left( cw\right) =\frac{3}{2}=1.5=\frac{\text{surface area}}{4}, 
\]

\[
\mathbb{E}\left( cw^{2}\right) =1+\frac{4}{\pi }=2.273239544735162..., 
\]

\[
\mathbb{E}\left( pw\right) =\frac{3}{2}\pi =4.712388980384689..., 
\]

\[
\mathbb{E}\left( pw^{2}\right) =8+6\pi \,_{3}F_{2}\left( -\tfrac{1}{2},%
\tfrac{1}{2},\tfrac{3}{2};1,2;1\right) =22.237117433439470..., 
\]%
\[
\mathbb{E}\left( cw\cdot pw\right) =2+\frac{16}{\pi }=7.092958178940650... 
\]%
which imply that the correlation between $cw$ and $pw$ is $0.915...$. \ We
shall present more details underlying these results later, especially that
involving the generalized hypergeometric function 
\[
_{3}F_{2}(a_{1},a_{2},a_{3};b_{1},b_{2};z)=\frac{\Gamma (b_{1})\Gamma (b_{2})%
}{\Gamma (a_{1})\Gamma (a_{2})\Gamma (a_{3})}\dsum\limits_{k=0}^{\infty }%
\frac{\Gamma (a_{1}+k)\Gamma (a_{2}+k)\Gamma (a_{3}+k)}{\Gamma
(b_{1}+k)\Gamma (b_{2}+k)}\frac{z^{k}}{k!}. 
\]%
whose appearance is quite unexpected.

\subsection{3-Crosspolytope (Octahedron)}

\[
\mathbb{E}\left( cw\right) =\frac{\sqrt{3}}{2}=0.866025403784438...=\frac{%
\text{surface area}}{4}, 
\]

\[
\mathbb{E}\left( cw^{2}\right) =\frac{1}{2}+\frac{\sqrt{2}}{\pi }-\frac{1}{%
2\pi }\func{arcsec}(3)=0.754244882063249..., 
\]

\[
\mathbb{E}\left( pw\right) =3\func{arcsec}(3)=3.692878252022324..., 
\]

\[
\mathbb{E}\left( pw^{2}\right) =13.6639421274..., 
\]%
\[
\mathbb{E}\left( cw\cdot pw\right) =3.2074623048... 
\]%
which imply that the correlation between $cw$ and $pw$ is $\approx 0.878$.
The fact that $cw_{3\text{-crosspolytope}}$ behaves like $2\cdot cw_{3\text{%
-simplex}}$ will be discussed shortly. We also have \cite{BoH}%
\[
\mathbb{P}\left( \text{projection has 4 vertices}\right) =\frac{3}{\pi }%
\left( \pi -2\func{arcsec}\left( 3\right) \right) =0.6490406878..., 
\]%
\[
\mathbb{P}\left( \text{projection has 6 vertices}\right) =\frac{2}{\pi }%
\left( 3\func{arcsec}\left( 3\right) -\pi \right) =0.3509593121.... 
\]%
See \cite{My1, My2, Vs} for related probabilities governing planar cross
sections of the tetrahedron, cube and octahedron (rather than projections).

\section{Tables of Intrinsic Volumes}

Let $\square $ be a rectangular $4$-parallelepiped in $\mathbb{R}^{4}$ of
dimensions $z_{1}$, $z_{2}$, $z_{3}$, $z_{4}$. It is well-known that \cite%
{Sntl}%
\[
V_{4}(\square )=z_{1}z_{2}z_{3}z_{4}, 
\]%
\[
V_{3}(\square
)=z_{1}z_{2}z_{3}+z_{1}z_{2}z_{4}+z_{1}z_{3}z_{4}+z_{2}z_{3}z_{4}, 
\]%
\[
V_{2}(\square
)=z_{1}z_{2}+z_{1}z_{3}+z_{1}z_{4}+z_{2}z_{3}+z_{2}z_{4}+z_{3}z_{4}, 
\]%
\[
V_{1}(\square )=z_{1}+z_{2}+z_{3}+z_{4} 
\]%
are the elementary symmetric polynomials in four variables. In $\mathbb{R}%
^{n}$, there are $n$ such intrinsic volumes, corresponding to the $n$
elementary symmetric polynomials. Limiting approximation arguments enable us
to compute $V_{j}(C)$ for arbitrary convex $C$.

One motivation for our work is to generalize the following two tables \cite%
{SW} for $n=2$:%
\[
\begin{array}{ccc}
V_{2}=\text{area,} &  & 2V_{1}=\text{circumference}%
\end{array}%
\]%
and $n=3$:%
\[
\begin{array}{ccccc}
V_{3}=\text{volume,} &  & 2V_{2}=\text{surface area,} &  & \tfrac{1}{2}V_{1}=%
\text{mean width}%
\end{array}%
\]%
to $n=4$:%
\[
V_{4}=\text{hypervolume,} 
\]%
\[
2V_{3}=\text{hyper-surface area,} 
\]%
\[
\tfrac{1}{3}V_{2}=\text{mean chorowidth,} 
\]%
\[
\tfrac{4}{3}V_{1}=\text{mean periwidth}=\pi \left( \text{mean width}\right) 
\text{.} 
\]

Another motivation is to use formulas for $V_{2}$ in \cite{BeH, HZ} to
deduce that 
\[
\mathbb{E}\left( cw_{n\text{-simplex}}\right) =\frac{n(n+1)}{8\sqrt{\pi }}%
\dint\limits_{-\infty }^{\infty }e^{-3x^{2}}\left( \frac{1+\func{erf}(x)}{2}%
\right) ^{n-2}\,dx 
\]%
for the regular $n$-simplex in $\mathbb{R}^{n}$; 
\[
\mathbb{E}\left( cw_{n\text{-cube}}\right) =\frac{n}{2} 
\]%
for the $n$-cube in $\mathbb{R}^{n}$; and 
\[
\mathbb{E}\left( cw_{n\text{-crosspolytope}}\right) =\frac{n(n-2)}{\sqrt{\pi 
}}\dint\limits_{0}^{\infty }e^{-3x^{2}}\func{erf}(x)^{n-3}\,dx 
\]%
for the\ regular $n$-crosspolytope in $\mathbb{R}^{n}$. In contrast, $%
\mathbb{E}\left( pw_{n}\right) $ are obtained simply by forming the product $%
\pi \cdot \mathbb{E}\left( w_{n}\right) $, and associated $\mathbb{E}\left(
w_{n}\right) $ are tabulated in \cite{Fi1, Fi2}. \ No such general
expressions are available for higher moments.

\section{Four-Dimensional Results}

From the foregoing, we calculate mean chorowidth and mean periwidth via
integration:%
\[
\mathbb{E}\left( cw\right) =\dint\limits_{0}^{2\pi }\dint\limits_{0}^{\pi
}\dint\limits_{0}^{2\pi }\dint\limits_{0}^{\pi }\dint\limits_{0}^{\pi }\text{%
area}\left( \text{convex\_hull}\left( M_{4}T\right) \right) \frac{1}{8\pi
^{3}}\sin \varphi \sin ^{2}\psi \sin \lambda \,d\psi \,d\varphi \,d\theta
\,d\lambda \,d\kappa , 
\]%
\[
\mathbb{E}\left( pw\right) =\dint\limits_{0}^{2\pi }\dint\limits_{0}^{\pi
}\dint\limits_{0}^{2\pi }\dint\limits_{0}^{\pi }\dint\limits_{0}^{\pi }\text{%
circumference}\left( \text{convex\_hull}\left( M_{4}T\right) \right) \frac{1%
}{8\pi ^{3}}\sin \varphi \sin ^{2}\psi \sin \lambda \,d\psi \,d\varphi
\,d\theta \,d\lambda \,d\kappa 
\]%
and likewise for mean square chorowidth, mean square periwidth and joint
moment.

\subsection{4-Simplex}

\[
\mathbb{E}\left( cw\right) =\frac{5\sqrt{3}}{12\pi }\left( \pi -\func{arcsec}%
(4)\right) =0.418889720727840..., 
\]

\[
\mathbb{E}\left( cw^{2}\right) =0.176..., 
\]

\[
\mathbb{E}\left( pw\right) =\frac{10}{3\pi }\left( 2\pi -3\func{arcsec}%
(3)\right) =2.748401146360593..., 
\]

\[
\mathbb{E}\left( pw^{2}\right) =7.56..., 
\]%
\[
\mathbb{E}\left( cw\cdot pw\right) =1.15... 
\]%
which imply that the correlation between $cw$ and $pw$ is $\approx 0.1$.
Closed-form expressions for $\mathbb{E}\left( cw^{2}\right) $, $\mathbb{E}%
\left( pw^{2}\right) $ and $\mathbb{E}\left( cw\cdot pw\right) $ are not
known. \ We also have%
\[
\mathbb{P}\left( \text{projection has 3 vertices}\right) \approx 0.146, 
\]%
\[
\mathbb{P}\left( \text{projection has 4 vertices}\right) \approx 0.585, 
\]%
\[
\mathbb{P}\left( \text{projection has 5 vertices}\right) \approx 0.269. 
\]%
These values are consistent with a theorem in \cite{AS} that the expected
number of vertices should be%
\begin{eqnarray*}
20\sqrt{\frac{2}{\pi }}\dint\limits_{-\infty }^{\infty }e^{-2x^{2}}\left( 
\frac{1+\func{erf}(x)}{2}\right) ^{3}dx &=&10\left( 1-\frac{3}{2\pi }\func{%
arcsec}(3)\right) =4.122... \\
&\approx &3(0.146)+4(0.585)+5(0.269).
\end{eqnarray*}

\subsection{4-Cube}

\[
\mathbb{E}\left( cw\right) =2, 
\]

\[
\mathbb{E}\left( cw^{2}\right) =4.04..., 
\]

\[
\mathbb{E}\left( pw\right) =\frac{16}{3}=5.\overline{3}, 
\]

\[
\mathbb{E}\left( pw^{2}\right) =28.4..., 
\]%
\[
\mathbb{E}\left( cw\cdot pw\right) =10.7... 
\]%
which imply that the correlation between $cw$ and $pw$ is $\approx 0.9$.

\subsection{4-Crosspolytope}

\[
\mathbb{E}\left( cw\right) =\frac{4\sqrt{3}}{9}=0.769800358919501..., 
\]

\[
\mathbb{E}\left( cw^{2}\right) =0.598..., 
\]

\[
\mathbb{E}\left( pw\right) =\frac{16}{\pi }\left( \pi -2\func{arcsec}%
(3)\right) =3.461550335020567..., 
\]

\[
\mathbb{E}\left( pw^{2}\right) =12.0..., 
\]%
\[
\mathbb{E}\left( cw\cdot pw\right) =2.67... 
\]%
which imply that the correlation between $cw$ and $pw$ is $\approx 0.8$. We
also have%
\[
\mathbb{P}\left( \text{projection has 4 vertices}\right) \approx 0.463, 
\]%
\[
\mathbb{P}\left( \text{projection has 6 vertices}\right) \approx 0.478, 
\]%
\[
\mathbb{P}\left( \text{projection has 8 vertices}\right) \approx 0.059. 
\]%
These values are consistent with a theorem in \cite{BoH} that the expected
number of vertices should be%
\begin{eqnarray*}
48\sqrt{\frac{2}{\pi }}\dint\limits_{0}^{\infty }e^{-2x^{2}}\func{erf}%
(x)^{2}dx &=&24\left( 1-\frac{2}{\pi }\func{arcsec}(3)\right) =5.192... \\
&\approx &4(0.463)+6(0.478)+8(0.059).
\end{eqnarray*}

\section{Further Work}

After having written the preceding, we discovered \cite{VB0}, which provides
new insights in the $3$-dimensional case. Earlier papers in this line of
thought include \cite{VB1, VB2, VB3, VB4}. \ Consider the $3$-cube with
vertices%
\[
T=\left( 
\begin{array}{cccccccc}
\tfrac{1}{2} & -\tfrac{1}{2} & \tfrac{1}{2} & \tfrac{1}{2} & -\tfrac{1}{2} & 
-\tfrac{1}{2} & \tfrac{1}{2} & -\tfrac{1}{2} \\ 
\tfrac{1}{2} & \tfrac{1}{2} & -\tfrac{1}{2} & \tfrac{1}{2} & -\tfrac{1}{2} & 
\tfrac{1}{2} & -\tfrac{1}{2} & -\tfrac{1}{2} \\ 
\tfrac{1}{2} & \tfrac{1}{2} & \tfrac{1}{2} & -\tfrac{1}{2} & \tfrac{1}{2} & -%
\tfrac{1}{2} & -\tfrac{1}{2} & -\tfrac{1}{2}%
\end{array}%
\right) . 
\]%
A starting point here is the simultaneous system of equations 
\[
\left\{ 
\begin{array}{l}
cw=x+y+z \\ 
pw=2\left( \sqrt{1-x^{2}}+\sqrt{1-y^{2}}+\sqrt{1-z^{2}}\right)%
\end{array}%
\right. 
\]%
given a unit vector $U=(x,y,z)$ in the first octant. This makes possible,
for example, the derivation of a closed-form marginal density for $cw$
(although not a joint density). Let us focus on computing $\mathbb{E}\left(
pw^{2}\right) $ and $\mathbb{E}\left( cw\cdot pw\right) $. Contributing to $%
\mathbb{E}\left( pw^{2}\right) $ are three terms like%
\[
I=\dint\limits_{0}^{\pi /2}\,\dint\limits_{0}^{\pi /2}\left( 1-\cos
^{2}\theta \sin ^{2}\varphi \right) \sin \varphi \,d\theta \,d\varphi =\frac{%
\pi }{3} 
\]%
and six terms like%
\begin{eqnarray*}
J &=&\dint\limits_{0}^{\pi /2}\,\dint\limits_{0}^{\pi /2}\sqrt{1-\sin
^{2}\theta \sin ^{2}\varphi }\,\sqrt{1-\cos ^{2}\varphi }\sin \varphi
\,d\theta \,d\varphi \\
&=&\dint\limits_{0}^{\pi /2}E\left( \sin \varphi \right) \sin ^{2}\varphi
\,d\varphi =\frac{\pi ^{2}}{8}\,_{3}F_{2}\left( -\tfrac{1}{2},\tfrac{1}{2},%
\tfrac{3}{2};1,2;1\right)
\end{eqnarray*}%
where%
\[
E(\xi )=\dint\limits_{0}^{\pi /2}\sqrt{1-\xi ^{2}\sin (\theta )^{2}}%
\,d\theta =\dint\limits_{0}^{1}\sqrt{\dfrac{1-\xi ^{2}t^{2}}{1-t^{2}}}\,dt 
\]%
is the complete elliptic integral of the second kind. The final result is $%
32(3I+6J)/(4\pi )$. For $\mathbb{E}\left( cw\cdot pw\right) $, the
calculations are simpler, with $I=\pi /6$ and $J=2/3$.

An analogous system of equations for the regular $3$-simplex with vertices \ 
\[
T=\left( 
\begin{array}{cccc}
0 & \frac{\sqrt{3}}{3} & -\frac{\sqrt{3}}{6} & -\frac{\sqrt{3}}{6} \\ 
0 & 0 & \frac{1}{2} & -\frac{1}{2} \\ 
\tfrac{\sqrt{6}}{4} & -\tfrac{\sqrt{6}}{12} & -\tfrac{\sqrt{6}}{12} & -%
\tfrac{\sqrt{6}}{12}%
\end{array}%
\right) 
\]%
is more difficult. Define constants%
\[
\begin{array}{ccc}
\gamma =\func{arccot}\left( 2\sqrt{2}\right) =0.339..., &  & \delta =\func{%
arccot}\left( \sqrt{2}\right) =0.615...%
\end{array}%
\]%
and functions%
\[
\begin{array}{ccc}
\alpha (\varphi )=\left\{ 
\begin{array}{lll}
\func{arcsec}\left( 2\sqrt{2}\tan \varphi \right) &  & \text{if }\gamma \leq
\varphi \leq \pi /2, \\ 
0 &  & \text{if }0\leq \varphi \leq \gamma ,%
\end{array}%
\right. &  & \beta (\varphi )=\dfrac{2\pi }{3}-\alpha (\varphi ).%
\end{array}%
\]%
Given a unit vector $U=(x,y,z)$ in the first dodecant (one of twelve regions
in $3$-space), we have%
\[
cw=\left\{ 
\begin{array}{lll}
\frac{1}{12}\left( \sqrt{6}x+3\sqrt{2}y+\sqrt{3}z\right) &  & \text{if }%
\delta \leq \varphi \leq \frac{\pi }{2}\text{ and }\beta (\varphi )\leq
\theta \leq \frac{\pi }{3}, \\ 
\frac{\sqrt{3}}{6}\left( \sqrt{2}x+z\right) &  & \text{if }\gamma \leq
\varphi \leq \frac{\pi }{2}\text{ and }0\leq \theta \leq \min \{\alpha
(\varphi ),\beta (\varphi )\}, \\ 
\frac{\sqrt{3}}{4}z &  & \text{if }0\leq \varphi \leq \delta \text{ and }%
\alpha (\varphi )\leq \theta \leq \frac{\pi }{3}%
\end{array}%
\right. 
\]%
and the three respective expressions for $pw$ are%
\[
\tfrac{1}{6}\left( 3\sqrt{4-\left( \sqrt{3}x-y\right) ^{2}}+\sqrt{3}\sqrt{%
12-\left( x-\sqrt{3}y+2\sqrt{2}z\right) ^{2}}+2\sqrt{3}\sqrt{3-\left( x-%
\sqrt{2}z\right) ^{2}}\right) , 
\]%
\begin{eqnarray*}
&&\tfrac{\sqrt{3}}{6}\left( \sqrt{3}\sqrt{4-\left( \sqrt{3}x-y\right) ^{2}}+%
\sqrt{3}\sqrt{4-\left( \sqrt{3}x+y\right) ^{2}}+\right. \\
&&\;\;\;\;\;\;\;\left. \sqrt{12-\left( x-\sqrt{3}y+2\sqrt{2}z\right) ^{2}}+%
\sqrt{12-\left( x+\sqrt{3}y+2\sqrt{2}z\right) ^{2}}\right) ,
\end{eqnarray*}%
\begin{eqnarray*}
&&\; \\
&&\tfrac{1}{2}\left( \sqrt{4-\left( \sqrt{3}x-y\right) ^{2}}+2\sqrt{1-y^{2}}+%
\sqrt{4-\left( \sqrt{3}x+y\right) ^{2}}\right) .
\end{eqnarray*}

An analogous system of equations for the regular $3$-crosspolytope with
vertices \ 
\[
T=\left( 
\begin{array}{cccccc}
\frac{1}{\sqrt{2}} & -\frac{1}{\sqrt{2}} & 0 & 0 & 0 & 0 \\ 
0 & 0 & \frac{1}{\sqrt{2}} & -\frac{1}{\sqrt{2}} & 0 & 0 \\ 
0 & 0 & 0 & 0 & \frac{1}{\sqrt{2}} & -\frac{1}{\sqrt{2}}%
\end{array}%
\right) 
\]%
is similar. Define functions%
\[
\begin{array}{ccc}
\alpha (\varphi )=\left\{ 
\begin{array}{lll}
\pi /4-\func{arcsec}\left( \sqrt{2}\tan \varphi \right) &  & \text{if }%
\delta \leq \varphi \leq \pi /2, \\ 
\pi /4 &  & \text{if }0\leq \varphi \leq \delta ,%
\end{array}%
\right. &  & \beta (\varphi )=-\alpha (\varphi )%
\end{array}%
\]%
where the constant $\delta $ is the same as before. Given a unit vector $%
U=(x,y,z)$ in the first hexadecant (one of sixteen regions in $3$-space), we
have%
\[
cw=\left\{ 
\begin{array}{lll}
\frac{1}{2}\left( x+y+z\right) &  & \text{if }\delta \leq \varphi \leq \frac{%
\pi }{2}\text{ and }\max \{\alpha (\varphi ),\beta (\varphi )\}\leq \theta
\leq \frac{\pi }{4}, \\ 
x &  & \text{if }\frac{\pi }{4}\leq \varphi \leq \frac{\pi }{2}\text{ and }%
0\leq \theta \leq \beta (\varphi ), \\ 
z &  & \text{if }0\leq \varphi \leq \frac{\pi }{4}\text{ and }0\leq \theta
\leq \alpha (\varphi )%
\end{array}%
\right. 
\]%
and the three respective expressions for $pw$ are%
\[
\sqrt{2}\left( \sqrt{2-(x+y)^{2}}+\sqrt{2-(x+z)^{2}}+\sqrt{2-(y+z)^{2}}%
\right) , 
\]%
\[
\sqrt{2}\left( \sqrt{2-(y-z)^{2}}+\sqrt{2-(y+z)^{2}}\right) ,\medskip 
\]%
\[
\sqrt{2}\left( \sqrt{2-(x-y)^{2}}+\sqrt{2-(x+y)^{2}}\right) . 
\]%
A rigorous proof that $cw_{3\text{-crosspolytope}}$ and $2\cdot cw_{3\text{%
-simplex}}$ are identically distributed remains open.

Let us return finally to (ordinary) width $w$ and questions left unanswered
in \cite{Fi3}. The\ $2$-cube (square) and regular $2$-simplex (equilateral
triangle) have vertices%
\[
\begin{array}{ccc}
\left( 
\begin{array}{cccc}
\tfrac{1}{2} & -\tfrac{1}{2} & \tfrac{1}{2} & -\tfrac{1}{2} \\ 
\tfrac{1}{2} & \tfrac{1}{2} & -\tfrac{1}{2} & -\tfrac{1}{2}%
\end{array}%
\right) , &  & \left( 
\begin{array}{ccc}
0 & \tfrac{1}{2} & -\tfrac{1}{2} \\ 
\tfrac{\sqrt{3}}{3} & -\tfrac{\sqrt{3}}{6} & -\tfrac{\sqrt{3}}{6}%
\end{array}%
\right)%
\end{array}%
\]%
respectively. Given a unit vector $(x,y)$ in the first quadrant, $w=x+y$ for
the square. Given a unit vector $(x,y)$ in the first sextant, $w=\frac{\sqrt{%
3}}{2}x+\frac{1}{2}y$ for the triangle. The corresponding width densities are%
\[
\begin{array}{ccc}
\left\{ 
\begin{array}{lll}
\dfrac{4}{\pi }\dfrac{1}{\sqrt{2-w^{2}}} &  & \text{if }1\leq w<\sqrt{2}, \\ 
0 &  & \text{otherwise,}%
\end{array}%
\right. &  & \left\{ 
\begin{array}{lll}
\dfrac{6}{\pi }\dfrac{1}{\sqrt{1-w^{2}}} &  & \text{if }\dfrac{\sqrt{3}}{2}%
\leq w<1, \\ 
0 &  & \text{otherwise}%
\end{array}%
\right.%
\end{array}%
\]%
respectively, via simple argument.

\section{Acknowledgements}

I am grateful to Rolf Schneider, Richard Vitale, Glenn Vickers and Daniel
Klain for their helpful correspondence. Much more relevant material can be
found at \cite{Fi4}, including experimental computer runs that aided
theoretical discussion here.

\bigskip

\end{document}